\documentclass[11pt]{amsart}
\usepackage{amsmath}
\usepackage{amsthm}
\usepackage{amsfonts}
\newtheorem{teo}{\bf Theorem}[section]

\newtheorem{cor}[teo]{\bf Corollary}

\newtheorem{preg}[teo]{\bf Question}

\newcommand{\Z}{\Cal{Z}}
\newcommand{\T}{\Cal{T}}
\newcommand{\F}{\Cal{F}}

\font\de=cmssi10
\begin{document}

\title[]{Dynamics in the isotopy class of a pseudo-Anosov map}
\thispagestyle{empty}
\begin{abstract}
Despite its homotopical stability, new relevant dynamics appear in the
isotopy class of a pseudo-Anosov homeomorphism. We study these new
dynamics by identifying homotopically equivalent orbits, obtaining
a more complete description of the topology of the corresponding quotient
spaces, and their stable and unstable sets. In particular, we get
some insight on how new periodic points appear, among other
corollaries. A list of further questions and problems is added at
the end of the paper.
\end{abstract}

\author{ Federico Rodriguez Hertz}
\address{CC 30, IMERL - Facultad  de Ingenier\'{\i}a\\Universidad de la
Rep\'{u}blica\\Montevideo, Uruguay}
\email{frhertz@fing.edu.uy}

\author{ Jana Rodriguez Hertz}
\address{CC 30, IMERL - Facultad  de Ingenier\'{\i}a\\Universidad de la
Rep\'{u}blica\\Montevideo, Uruguay}
\email{jana@fing.edu.uy}

\author{ Ra\'{u}l Ures}
\address{CC 30, IMERL - Facultad  de Ingenier\'{\i}a\\Universidad de la
Rep\'{u}blica\\Montevideo, Uruguay}
\email{ures@fing.edu.uy}
\thanks{Authors were supported in part by a CONICYT-Clemente
Estable grant}

\maketitle
%

\renewcommand{\theenumi}{\arabic{section}.\arabic{enumi}}
\renewcommand{\theequation}{\arabic{section}.\arabic{equation}}
%
\newcommand{\noteq}{\neq}
\newcommand{\eps}{\varepsilon}
\newcommand{\aster}{\ast}
\newcommand{\zeda}{\zeta}
\newcommand{\etha}{\eta}
\newcommand{\xsi}{\xi}
\newcommand{\uperline}{\bar}
\newcommand{\Cal}{\mathcal}
\newcommand{\disp}[1]{\displaystyle{\mathstrut#1}}
\newcommand{\fra}[2]{\displaystyle\frac{\mathstrut#1}{\mathstrut#2}}
\renewcommand{\~}{\widetilde}
%


\section{\bf Introduction}
In this paper we consider a pseudo-Anosov homeomorphism $f:M^2\to
M^2$ on a closed surface, and homeomorphisms $g:M^2\to M^2$ in its
isotopy class. The reader not familiar with these concepts may consult \cite{flp}.
Following \cite{h}, we shall use Katok's {\de global shadowing}
to compare orbits of $f$ with orbits of $g$. According to this
definition the $f$-orbit of $x$ is globally shadowed by the
$g$-orbit of $y$ if $\sup_{n\in\mathbb Z}D(F^n(\~x), G^n(\~y))<\infty$,
for some adequate lifts $F,G,\~x,\~y$ of $f,g,x,y$, respectively,
to the universal cover ($D$ is any equivariant metric). For simplicity we denote $(f,x)\sim(g,y)$.\par
This equivalence relation permitted to state a kind of homotopical
stability of $f$ in its isotopy class, in the sense that for all $g\simeq f$
each $f$-orbit $(f,x)$ is globally shadowed by (at least) one $g$-orbit $(g,y_x)$
\cite{l1,h,f}. That is, the dynamics of $f$ is repeated all over
its isotopy class. However, unlike the Anosov case, new dynamics
may appear for $g$, which are not globally shadowed by $f$.\par
In this work we study these new dynamics using the following
approach. Given a fixed $g$, we identify all $(g,x)\sim(g,x')$,
obtaining a quotient space $M_L(g)$. The map $g$ naturally induces a homeomorphism $g_L$ on $M_L(g)$.
Observe that, in general $M_L(g)$ is no longer a surface, though $g_L$ is still expansive.\par
 The core of this paper is based on an observation,
due to Fathi \cite{f}, of the existence of a ``universal torus" ${\mathcal Z}$ and a hyperbolic extension $\Phi$ of $f$ in ${\mathcal Z}$, for
which $M$ is the smallest non trivial invariant sub-continuum.
The universal cover of ${\mathcal Z}$ is the product of
two trees obtained from the stable and unstable foliations,
respectively.\par
Our first theorem states that all $M_L(g)$ are naturally embedded
in ${\mathcal Z}$, and that $g_L$ is nothing but $\Phi$ restricted to $M_L(g)$
\begin{teo}\label{teo1}
Let $f:M^2\to M^2$ be a pseudo-Anosov homeomorphism on a closed
surface. Then there exists ${\mathcal Z}$ and a homeomorphism $\Phi:{\mathcal Z}\to{\mathcal
Z}$, such that:
\begin{enumerate}
\item\label{teo1.1} $\Phi:{\mathcal Z}\to{\mathcal Z}$ is a {\de metrically split
 hyperbolic homeomorphism}.
\item \label{teo1.2}$\Phi$ has the pseudo-orbit tracing property on
$\Z$. Additionally, the pseudo orbits of any lift of
$\Phi$ are uniquely shadowed in the universal cover of $\Z$.
\item \label{teo1.3}For all $g:M^2\to M^2$ homotopic to $f$, there exists a
$\Phi$-invariant subset $L(g)\subset\Z$ on which $g_L$ is
conjugate to $\Phi|_{L(g)}$
\end{enumerate}
\end{teo}
\subsection*{Remark}
A {\de hyperbolic homeomorphism} is to be understood in the sense of Bowen, that is, one having hyperbolic
canonical coordinates \cite{b2}. The definition of {\de metrically splitting} is, {\it mutatis
mutandis}, as in Franks (\cite{fr}, p.67).\newline\par

Theorem \ref{teo1} provides a unified scenario for all $M_L(g)$ and $g_L$, endowed with a hyperbolic
structure. From now on, we shall write $L(g)$ and $\Phi_g=\Phi|_{L(g)}$ instead of $M_L(g)$ and $g_L$.
We find that:
\begin{teo}\label{cont.coc} The application $g\mapsto L(g)\in 2^\Z$ is
continuous in the isotopy class of $f$, $2^\Z$ endowed with the
Hausdorff topology. The set $\{L(g):g\simeq f\}\subset
2^\Z$ is arcwise connected
\end{teo}
As a Corollary we obtain
\begin{cor}\label{top.est}
$\Phi_g$ is never topologically stable. If $g$ and $g'$ are
semiconjugate by a semiconjugacy homotopic to the identity then
$L(g)=L(g')$.
\end{cor}
\par
Corollary \ref{top.est} shows that each $\Phi_g$ has a distinct
dynamics.
 Let us point out that always when a pair $g$, $g'$ are
 semiconjugate, a semiconjugacy can be chosen to be homotopic to a
 symmetry of $f$ (see comments in \S \ref{sec.preg}). Although previous statements give some light to the new dynamics
appearing in the isotopy class of $f$, they mainly open many
questions and problems, in \S\ref{sec.preg} we state some of
them.\par
It is expectable that the study, and specially the local
analysis of the sets $L(g)$ will bring information about the
dynamics of the $\Phi_g$, and are an interesting topic in its own.
Following result describes the local properties of $L(g)$, new
periodic points, and local stable and unstable sets.

 \begin{teo}\label{teo.local}Let $g:M^2\to M^2$ be a homeomorphism on
 a closed surface, $g\simeq f$. Then:
 \begin{enumerate}
 \item \label{teo4.1} The topological dimension of $L(g)$ is $2$.
 \item \label{teo4.2} There exists $\eps>0$ such that $CW^s_{\varepsilon}(x)$
 and $CW^u_\eps(x)$ are nontrivial dendrites for all $x\in L(g)$,
 with diameters bounded from below.
 \item \label{teo4.3}If $p\in{\mathcal P}er(\Phi_g)\setminus{\mathcal
 P}er(\Phi_f)$ then $p$ is an endpoint of {\de both}
 $CW^s_{\varepsilon}(p)$ and $CW^u_\eps(p)$.
 \end{enumerate}
 \end{teo}
The {\de $\eps$-stable component of $x$}, $CW^s_\eps (x)$, denotes the
component of $x$ in $W^s_\eps(x)$, the $\eps$-local stable set of
$x$. A dendrite is a uniquely arcwise connected curve.\newline\par
 Item \ref{teo4.2} improves \cite{cs} where, using delicate arguments of plane topology,
it was shown arcwise connectedness of $CW^\sigma_\eps(x)$ for $g$ $C^0$-close to
an expansive homeomorphism. It also extends the result over the isotopy class of
$f$.
\subsection*{Acknowledgements} This paper is related to questions and problems posed by J.
Lewowicz. We want to express our acknowledgement to him for
bringing the problem addressed in Theorem \ref{teo.local} to our attention.
%
\section{\bf Global Properties}\label{seccion.global}
%
From now on, $f:M\rightarrow M$ will be a pseudo-Anosov homeomorphism and
$F:\~{M}\rightarrow \~{M}$ any lift of $f$ to the universal covering of $M$.
 Associated to $F$ there are two transversal invariant foliations with singularities,
 $\Cal F^{u}$ and $\Cal F^{s}$, with transverse measures contracted and expanded by $F$
 with rates $\lambda>1$ and $\lambda^{-1}$, respectively.\par
We briefly describe some remarkable properties of the construction used by A. Fathi to embed
a pseudo-Anosov homeomorphism in a hyperbolic dynamics. The reader
interested in the details should consult \cite{f}.
\begin{enumerate}
\item \label{item.1} The set of $\sigma$-leaves $\T^\sigma=\~M/_{\~\F^\sigma}$ is a
non locally compact, not complete tree, when the distance of two points in $\T^\sigma$ is defined as the
minimum over the transverse measures of arcs joining their associated
leaves ($\sigma=s,u$)
\cite{ms}.\vspace*{0.5em}
\item The completion $\widehat{\T}^\sigma$ of $\T^\sigma$ is still a non
locally compact tree \cite{ms}, on which $F$ and the fundamental
group $\Gamma$ of $M$ naturally act, since the covering
transformations of $F$ are isometries of the metrics defined in
(\ref{item.1}).\vspace*{0.5em}
\item The metric associated to $\widehat{\T}^s$ is $\lambda$-expanded
under the action of $F$, while the metric of $\widehat{\T}^u$ is
$\lambda^{-1}$-contracted.\vspace*{0.5em}
\item The product action $\~\Phi$ on $\~\Z=\widehat{\T}^s\times\widehat{\T}^u$ is a {\de metrically split hyperbolic
homeomorphism}. $\Z=\~\Z/_{\Gamma}$ is a metric space having
$\~\Z$ as universal cover (the product metric is
$\Gamma$-equivariant). $\~\Phi$ trivially induces a homeomorphism $\Phi$ on $\Z$.\vspace*{0.5em}
\item \label{item2.5} The set $\~{L(f)}=\{(l^s,l^u)\in\~\Z:l^s\cap
l^u\ne\emptyset\}$ is homeomorphic to $\~M$. Its projection to $\Z$, $L(f)$, is
$\Phi$-invariant and $\Phi_f=\Phi|_{L(f)}$ is conjugate to $f$.
Moreover, $L(f)$ is contained in all non trivial $\Phi$-invariant closed
connected subsets of $\Z$ \cite{f}.
\end{enumerate}

Bowen's proof of the shadowing lemma \cite{b}
works also in this setting, both for $\Phi$ and for $\~\Phi$, since $\Z$ and $\~\Z$ are complete.
Being $\~\Z$ metrically split, the fact that $\~\Phi$ is a
hyperbolic homeomorphism easily implies that every
$\~\Phi$-pseudo-orbit is uniquely shadowed.\par
Let us prove that each $g_L$ is a restriction of $\Phi$ to a set homeomorphic to $M_L(g)$:

In view of \ref{item2.5}, $\~{L(f)}$ is a universal cover of $M$
and $\~\Phi$ is a lift of $f$ when restricted to this set. Take
 a lift $G$ of $g$, equivariantly homotopic to
 $\~\Phi|_{\~{L(f)}}$. Then each $G$-orbit is a pseudo
 orbit of $\~\Phi$. Shadowing Lemma provides a continuous map $\~H:\~{L(f)}\to \~\Z$
that is a $\Gamma$-invariant semiconjugacy between $G$ and
$\~\Phi|_{Im(\~H)}$. This induces a semiconjugacy
$H:L(f)\to Im(H)\subset\Z$, with the property that
$x'\in H^{-1}(x)$ if and only if $(g,x)\sim(g, x')$. Hence
$Im(H)=L(g)$ is homeomorphic to $M_L(g)$ and $g_L$ is
conjugate to $\Phi|_{L(g)}$, ending proof of Theorem
\ref{teo1}.\par

It is easy to prove that semiconjugacies that send $M$ into $\Z$
vary continuously with respect to homeomorphisms $g$ of $M$, in
the $C^0$ topology. This implies that $\Cal L=\{L(g)\subset \Z;\:
g\,\,\mbox{homeomorphism of}\,\,M\}$ is an arcwise-connected set
in the Hausdorff topology. In particular, no $L(g)$ is an isolated set of
$\Phi$. This implies they have not the shadowing property, which
in this context means they are not topologically stable.\par
For the proof of Corollary \ref{top.est}, take a continuous $h:M\to M$ homotopic
to the identity such that $h\circ g=g'\circ h$. Consider a lift
$\~h$ of $h$ to $\~{L(f)}$, and the $\Gamma$-invariant
semiconjugacies $\~H, \~H':\~{L(f)}\to \~\Z$ obtained in the proof
of Theorem \ref{teo1.3} for lifts $G$, $G'$ of $g$, $g'$
respectively. We claim that $\~H=\~H'\circ \~h$ for a suitable $\~h$. Indeed,
as $\~h$ can be chosen to be at a bounded distance of the identity, we easily
obtain that the distance between $\~\Phi^n\circ\~H(x)(=\~H\circ G^n(x))$
and $\~\Phi^n\circ \~H'\circ \~h(x)(=\~H'\circ \~h\circ G^n(x)))$ is
uniformly bounded over $n\in{\mathbb Z}$. We conclude $\~H(x)=\~H'\circ \~h(x)$
due to hyperbolicity of $\~\Phi$. The fact that $\~h$ is onto
implies $Im(H)=Im(H')$ and the Corollary follows.

\section{Local Properties}

This section is devoted to proving Theorem \ref{teo.local}:
\noindent\begin{enumerate}

\item Item \ref{teo4.1} is in part a consequence of Theorem
\ref{teo1.3}, since $\Z$, being the product of two
trees, has topological dimension 2. On the other hand, $L(g)$
contains a copy of $M$ (item (\ref{item2.5}) of Section
\S\ref{seccion.global}, see also \cite{l1,h}), hence its topological dimension is 2.

\item Item \ref{teo4.2} follows trivially from the fact that $CW^s_{\varepsilon}(x)$ are compact connected
subsets of a tree. These sets have uniform size due to \cite{l1} ($L(g)$ are locally connected).

\item Let $p$ be a periodic point of $\Phi$. Take $\~p=(l^s,l^u)$ a lift of $p$.  If $l^s\in \T^s$,
there exist a lift $\~\Phi$ of $\Phi$ and a covering transformation $\tau$ such
that $\tau \~\Phi^n(\~p)=\~p$, implying $l_s$ is periodic for $f:M\rightarrow M$. It follows that $l^s$
contains an n-periodic $q$ of $f$. There exists a lift $\~q$ of $q$ to $\~{L(f)}$
that is fixed by $\tau\~ \Phi^n$. Since $\tau \~\Phi^n_{\Z}$ has a unique fixed point we obtain that $p=q\in
M$. This shows a $\Phi$-periodic point is in $L(f)$ if and only if its corresponding leaves $l^s_p$ {\de and} $l^u_p$ are,
respectively, in the trees $\T^s$ and $\T^u$. \par
Now, points in the completion $\hat\T^\sigma$ that are not in $\T^\sigma$, are endpoints of the tree
\cite{f}, so we arrive to the desired conclusion.

\end{enumerate}
\section{Questions \& Comments}\label{sec.preg}
The first set of questions concerns the relation between
topological entropy, periodic points, and the sets $L(g)$
\begin{preg}
Is the application $g\mapsto h_{top}(\Phi_g)$ continuous? Is it
strictly increasing with respect to inclusions of $L(g)$'s?
\end{preg}
\begin{preg}
Does Bowen's formula apply for $h_{top}(\Phi_g)$?
\end{preg}
\begin{preg}\label{p.per}
Are periodic points dense in $L(g)$?
\end{preg}
Let us point out that, in case $g$ is Axiom A with the strong transversality
condition, it is known from \cite{lu} that question \ref{p.per} has an affirmative
answer.\par
Observe that $L(g)=L(g')$ does not necessarily imply $g$ is
semiconjugate to $g'$. However; an equivalence relation can be
defined for $g$ in the isotopy class of $f$, namely: $g\approx g'$
iff $L(g)=L(g')$. Is this relation dynamically relevant in some
sense? For example:
\begin{preg}
Is every $L(g)$ realized by some Axiom A? That is, does every
$\approx$-class contain an Axiom A?
\end{preg}
\begin{preg}
Is there a minimal representative (e.g. in the sense of its number of periodic points, topological entropy, etc) in each
$\approx$-class?
\end{preg}
\begin{preg} Does any $\Phi$-invariant subcontinuum of $\Z$
correspond to a $L(g)$, for some $g$ in the isotopy class of $f$?
\end{preg}
We would like also to observe that this question is still
unsolved:
\begin{preg}
Is $\dim_p L(g)=2$ at {\em all} points $p\in L(g)$? Does $L(g)$
have a special topological structure, such as Cantor manifold,
etc.?
\end{preg}
Finally, though it is not in the scope of our paper, it is
interesting to note that the existence of two semiconjugate
homeomorphisms in the isotopy class of a pseudo-Anosov map $f$, whose
semiconjugacy cannot be chosen homotopic to the identity indicates
the presence of symmetries of $f$.


\end{document}